\documentclass[titlepage,12pt]{article}
\usepackage{amssymb}
\usepackage{amsfonts}
\usepackage{amsmath}
\begin{document}

{\Large \bf A Problem in Categories} \\ \\

{\bf Elem\'{e}r E Rosinger} \\
Department of Mathematics \\
and Applied Mathematics \\
University of Pretoria \\
Pretoria \\
0002 South Africa \\
eerosinger@hotmail.com \\ \\

{\bf Abstract} \\

The problem is posed to find out for arbitrary nonvoid sets $X$ which are all the mappings
$T : X \longrightarrow X$ that can be defined and each separately identified through means of
categories alone. As argued, this problem may have a certain foundational relevance. \\ \\

{\bf 1. The Problem} \\

Find out which mappings $T : X \longrightarrow X$, with arbitrary nonvoid sets $X$, can be
defined and identified each separately by means of categories only. \\ \\

{\bf 2. Examples} \\

1) If $T = id_X$ is the {\it identity} mapping of $X$, then for every two mappings $f, g : X
\longrightarrow Y$, where $Y$ is an arbitrary set, one has the {\it cancellation} property \\

(1)~~~ $ f \circ T = g \circ T ~~~\Longrightarrow~~~ f = g $ \\

while for every two mappings $f, g : Y \longrightarrow X$, where $Y$ is an arbitrary set, one
has the dual {\it cancellation} property \\

(2)~~~ $ T \circ f = T \circ g ~~~\Longrightarrow~~~ f = g $ \\

In terms of categories, the identity mapping $T = id_X$ has of course the axiomatic
property \\

(3)~~~ $ f \circ T = f,~~~ T \circ g = g,~~~
                   f : X \longrightarrow Y,~~ g : Y \longrightarrow X $ \\

from which (1) and (2) result immediately. However, the question remains to what extent is the
identity mapping $T = id_X$ characterized by (3), or for that matter, (1) and (2), in terms of
categories only. \\

2) If $T$ is a {\it constant} mapping, that is, for a certain $c \in X$, we have $T ( x ) =
c$, with $x \in X$, then for every two mappings $f, g : Y \longrightarrow X$, where $Y$ is an
arbitrary set, one has the {\it coequalizer} property \\

(4)~~~ $ T \circ f = T \circ g $ \\

We note however that, while (4) may happen to {\it define} the set of constant mappings $T :
X \longrightarrow X$ as a whole, it certainly does {\it not} in general {\it identify} them
individually as well. \\

Therefore, the Problem above has in fact two {\it subproblems} :

\begin{quote}

I)~~~ Define by means of categories the largest class of mappings $T : X \longrightarrow X$,
where $X$ is an arbitrary set. \\

II)~~ Identify individually by means of categories the largest class of mappings $T : X
\longrightarrow X$, where $X$ is an arbitrary set.

\end{quote}

{~} \\

{\bf 3. On the Relevance of the Problem} \\

As far as the author is concerned, he has not seen the above Problem formulated, let alone
solved anywhere in the literature. The relevance of the Problem, in case it has indeed not
been considered before, may be {\it foundational}, as argued in what follows. \\

Category Theory, as introduced in [1], and typically presented ever since in the respective
literature, starts from Set Theory which is assumed to be given, and then follows with the
definition of categories through certain axioms formulated in set theoretic terms. \\
However, from foundational point of view, this approach is not the only one which has been
considered in the literature, [3, pp. 235-250]. In particular, the position of sets, versus
categories, when seen in a foundational perspective, can be changed, with categories being
considered as given, and sets being introduced in terms of categories. \\

In that latter case, however, the question arises to what extent can one recover, purely in
terms of categories, the structural richness involved in each and every {\it specific}
set, as inherent in it, when considered with Set Theory ? \\

And obviously, for any given set $X$, one of the immediate and naturally associated structures
is that of the set $X^X$ of all mappings $T : X \longrightarrow X$. \\

In this way, the above Problem does indeed address the foundational issue of whether there
exists the possibility of recovering the specific structural richness of Set Theory, and
recovering it in terms of Category Theory alone. \\

\end{document}